\documentclass[a4paper,10pt,oneside]{article}

\def\TITLE{Probabilities are always axiomatizable}
\def\AUTHOR{Zal\'an Gyenis \\[2mm] Jagiellonian University}
\def\DATE{\today}
\def\ABSTRACT{In this paper we study the interaction between logic and probability. In particular, we show that the convex hull of evaluations of a broad class of logics is always effectively axiomatizable.
We define a Birkhoff-style calculus for probability axioms for which compactness, and finite completeness is proved. We give example for a logic for which probabilities are not finitely axiomatizable.}
\def\KEYWORDS{} 
\def\SUBJCLASS{}   

\newif\ifpdfrender		\pdfrendertrue		
\newif\ifdraft			\draftfalse			
\newif\ifshowkeys		\showkeysfalse		
\newif\ifhyperlinked	\hyperlinkedfalse	
\newif\ifswap			\swaptrue			

\def\enddefsymbol{$\square$}				
\def\endproofsymbol{$\blacksquare$}			

\usepackage[utf8]{inputenc}
\usepackage[T1]{fontenc}
\usepackage[english]{babel}

\usepackage{amsmath,amssymb, graphicx, amsopn, amsthm}
\usepackage{enumerate, framed}
\usepackage[explicit]{titlesec}
\usepackage{nicefrac}	

\usepackage[dvipsnames]{xcolor}

\usepackage{microtype} 
\microtypesetup{protrusion=true,expansion=true}

\usepackage[colorinlistoftodos]{todonotes}
\newcommand{\TODO}[1]{\todo[color=green!40]{#1}}

\usepackage[numbers]{natbib}
\setcitestyle{square}

\usepackage{type1cm}         
\usepackage{makeidx}         
\usepackage{graphicx}        
\usepackage{multicol}        
\usepackage[bottom]{footmisc}

\usepackage{setspace}
\ifhyperlinked
	\usepackage[hidelinks]{hyperref} 
	\hypersetup{
		colorlinks,
    	linkcolor={blue!50!black},
    	citecolor={blue!50!black},
    	urlcolor={blue!80!black}
		}
	\usepackage[doipre={DOI:~}]{uri} 
\fi

\newcommand{\MAKETITLE}
{
\title{\textsc{\textbf{\MakeLowercase{\TITLE}}}}
\author{\AUTHOR}
\date{\DATE}
\maketitle
\begin{abstract}
	\noindent \ABSTRACT

	\ifx\KEYWORDS\empty\else
	  \vspace{5mm}
	  \noindent {\bf Keywords:} \KEYWORDS.
	\fi
	\ifx\SUBJCLASS\empty\else
	
	  \vspace{1mm}
	  \noindent {\bf Subject classification:} \SUBJCLASS.
	\fi
\end{abstract}

\vspace{5mm} \normalsize	
}

\ifpdfrender
	\usepackage{pdfrender}
	\pdfrender{StrokeColor=black,TextRenderingMode=2,LineWidth=0.15pt}
\fi

\ifdraft
	\usepackage{scrtime} 
	\usepackage[draft]{prelim2e}

\fi

\ifshowkeys
	\usepackage[notref, notcite]{showkeys}
	
\fi

\allowdisplaybreaks				

\usepackage[overload]{textcase} 	

\newcommand{\gyzFormatSec}[1]{\textbf{\textsc{\MakeLowercase{#1}}}}

\titleformat{\section}
{\normalfont}{\Large\fontfamily{yvtj}\selectfont\thesection}{1em}{\Large\gyzFormatSec{#1}}

\titleformat{\subsection}
{\normalfont}{\large\fontfamily{yvtj}\selectfont\thesubsection}{1em}{\large\gyzFormatSec{#1}}

\makeatletter
\let\osection\section
\def\section{\@tempskipa\lastskip\removelastskip
\penalty-20
\vskip\@tempskipa\osection}
\let\osubsection\subsection
\def\subsection{\@tempskipa\lastskip\removelastskip
\penalty-20
\vskip\@tempskipa\osubsection}
\makeatother

\ifswap
\newtheoremstyle{mythmstyle}
  {\topsep}
  {\topsep}
  {\normalfont}
  {}
  {}
  {{\bf .}}
  {.7em}
  {{{\bfseries\thmnumber{#2}}~{\bfseries\thmname{#1}}{\normalfont{\thmnote{ (#3)}}}}} 
\else
\newtheoremstyle{mythmstyle}
  {\topsep}
  {\topsep}
  {\normalfont}
  {}
  {}
  {\bfseries.}
  {.7em}
 {{{\bfseries\thmname{#1}}~{\bfseries\thmnumber{#2}}{\normalfont{\thmnote{ (#3)}}}}} 
\fi

\theoremstyle{mythmstyle} 
\newtheorem{theorem}{Theorem}[section]
\newtheorem{lemma}[theorem]{Lemma}

\newtheorem{proposition}[theorem]{Proposition}
\newtheorem{corollary}[theorem]{Corollary}

\newtheorem{observation}[theorem]{Observation}

\newtheorem{defi/}[theorem]{Definition}
\newenvironment{definition}
  {\begin{defi/}}
  {\hfill\enddefsymbol
	\end{defi/}
  
  	\bigskip}

\newtheoremstyle{mynotestyle}
  {\topsep}
  {\topsep}
  {\normalfont}
  {}
  {}
  {\bfseries.}
  {.7em}
  {{{\bfseries\thmnumber{#2}}{\normalfont{\thmnote{ (#3)}}}}}

\theoremstyle{mynotestyle}

\renewenvironment{proof}[1][\unskip]{%
\par
\noindent
\textbf{Proof #1.}
\noindent}
{\hfill\endproofsymbol

\bigskip}

\let\<=\langle
\let\>=\rangle

\def\[#1\]{\begin{align}#1\end{align}}

\def\RED{\textcolor{red}}
\def\BLUE{\textcolor{blue}}

\def\changeA{\marginpar[\hfill $\downarrow$Change]{Change$\downarrow$}}
\def\changeZ{\marginpar[\hfill $\uparrow$Change]{Change$\uparrow$}}

\def\land{\wedge}       
\def\lor{\vee} 
\def\Land{\bigwedge}    
\def\Lor{\bigvee}
\def\lnot{\neg}

\let\models=\vDash

\def\proves{\vdash}

\def\defeq{\overset{\text{def}}{=}}

\DeclareMathOperator{\Hom}{Hom}

\def\HSP{\mathbf{HSP}}

\def\gA{\mathfrak{A}}

\def\gF{\mathfrak{F}}

\def\gL{\mathfrak{L}}

\def\gT{\mathfrak{T}}

\def\cL{\mathcal{L}}

\def\cU{\mathcal{U}}

\def\FEQ{\mathsf{FEQ}}

\DeclareMathOperator{\Ineq}{Ineq}
\DeclareMathOperator{\Id}{Id}
\DeclareMathOperator{\CH}{CH}
\def\TWO{\mathbf{2}}

\def\TRUE{1}
\def\FALSE{0}

\def\MIDDLEi{\nicefrac{1}{2}}
\def\restr{\!\upharpoonright\!}

\def\mathL{\text{\it{\L}}}

\let\phi=\varphi
\let\theta=\vartheta

\usepackage[left=2.6cm,right=2.6cm]{geometry}

\begin{document}
	\MAKETITLE

\section{The theme}

Belief is often represented by a function assigning real numbers from the interval
$[0,1]$ to propositions. If complete and perfect information were available, beliefs could
simply match classical truth-values, taking value $1$ for truths and $0$ for falsehoods. In practice, however, information is limited or uncertain. By classical theorems of Ramsey and de Finetti, rational belief functions—that is, belief functions avoiding Dutch books—are exactly those functions that satisfy the Kolmogorov axioms of probability. It is well known that probabilities can be described as convex combinations of classical truth-valuations. Work by Paris \cite{Paris2005}, Williams \cite{Williams2016},
and Bradley \cite{Bradley2016} has shown that this characterization also holds for various kinds of nonclassical logic. Axiomatizing probabilities and reasoning with belief functions is the topic of several recent papers: \cite{Bilkova2023BelnapDunn} for Belnap--Dunn logic, \cite{Gerla2000MVAlgebras} for MV-algebras, and \cite{GSGyW2022a,GSGyW2022} for symmetric logic.

In this paper we show that probability functions over propositional-like logics are always axiomatizable. Finding such an axiomatization is effective, and for logics with finitely many variables we can always obtain a finite axiom system. We also give an example of a logic for which there is no finite set of axioms describing probabilities (thus, in general, we cannot expect ``compactness''). We build a Birkhoff-style calculus for probability axioms, for which we prove compactness and completeness.

The structure of the paper is as follows. In the next section we briefly summarize earlier results that are important for us, and then in Section \ref{sec:genpat} we discuss the scope of these results.
Section \ref{sec:new} contains the main theorems of the paper. We define an equational calculus for belief functions; then in Section \ref{sub:sep} we introduce conditions that provide an equational axiomatization of probabilities, generalizing many earlier known results. In Subsection \ref{sub:always} we prove that such axiomatizations always exist. Finally, Subsection \ref{subsec:nonfin} presents an example of a logic in which probabilities can be axiomatized, but not finitely.

\section{Summary of previous results}\label{sec:prelim}
By the classical theorem of Ramsey and De Finetti (see e.g. \cite{DeFinetti1974}) a belief function is rational, in the sense of avoiding any Dutch book, if and only if it satisfies the axioms of probability. Let us unfold the details. Take a set $P$ of propositional letters, and using the standard logical connectives $Cn = \{\lor, \land, \lnot, \to\}$ generate the formula algebra $\gF = \gF(P, Cn)$ of propositional logic. Let $\TWO$ be the two-element Boolean algebra with universe $\{0,1\}$. Evaluations of propositional logic are the homomorphisms $v:\gF\to\TWO$, and the semantic consequence relation $\models$ is defined as $\varphi\models\psi$ if and only if $v(\varphi)=1$ implies $v(\psi)=1$ for all homomorphism $v\in \Hom(\gF, \TWO)$. A belief function is any function $B:\gF\to[0,1]$. 
The standard Dutch book argument goes as follows.
Suppose that we identify one's degree of belief $B(\varphi)$ in a certain formula $\varphi$ with the ``betting quotient'' at which the person, call it bettor, is ready to bet that $\varphi$ is true. The bookie asks the bettor to give his betting quotient $B(\varphi)$. Then the bookie sets his stake $s(\varphi)$. The bettor pays the bookie the sum $B(\varphi)\cdot s(\varphi)$ in exchange for the right to receive the amount $s(\varphi)$ from the bookie if $\varphi$ is true. Otherwise, the bettor gets nothing, that is, looses the amount $B(\varphi)\cdot s(\varphi)$. The bettor's betting quotients $B(\varphi_1)$, $\ldots$, $B(\varphi_n)$
are called coherent, if the bookie cannot chose the stakes $s(\varphi_1)$, $\ldots$, $s(\varphi_n)$ in such a way that the bettor wins irrespective of which of the formulas $\varphi_1$, $\ldots$, $\varphi_n$ are true. The bookie makes a Dutch book against the bettor, if the bookie can chose the stakes so that he wins whichever formula becomes true. Then, the Ramsey and De Finetti theorem can be formulated as follows.

\begin{theorem}[Cf. \cite{DeFinetti1974,Paris2005}]
	There is no Dutch book against the bettor if and only if $B$ satisfies the axioms of probability. More precisely, the following are equivalent.
	\begin{enumerate}[(A)]\itemsep-3pt
		\item {\bf No Dutch book}:  There are no $\theta_1,\ldots,\theta_n\in\gF$ and $s_1,\ldots, s_n\in\mathbb{R}$
		such that for all evaluations $v\in\Hom(\gF, \TWO)$ we have
		\[	\sum_{i=1}^{n}s_i\cdot(v(\theta_i)-B(\theta_i)) \ \ < \ \ 0  \]
		Here $s_i$ is the price associated with the bet for $\theta_i$, $s_i\cdot B(\theta_i)$ is the price of that bet, and
		$s_i\cdot v(\theta_i)$ is the payout if $\theta_i$ has truth value $v(\theta_i)$  ($1$ means true, and $0$ means false).		
		\item {\bf Probability}: $B$ is a probability function if for all $\varphi, \psi\in\gF$:
		\begin{itemize}
			\item[(P1)] If $\models\phi$, then $B(\phi)=1$
			\item[] If $\models\lnot\phi$, then $B(\phi)=0$
			\item[(P2)] If $\phi\models\psi$, then $B(\phi)\leq B(\psi)$
			\item[(P3)] $B(\phi\lor\psi) = B(\phi)+B(\psi)-B(\phi\land\psi)$
			\qed
		\end{itemize}
	\end{enumerate}
\end{theorem}

\noindent In the larger context of norms of rationality, one then typically applies the classical Dutch Book justification of identifying rational belief functions with probabilities (even in nonclassical settings, see later). This argument and its justification are not topics of the present paper; we refer the Reader to \cite{Paris2005,Williams2016,Bradley2016}. 

Evaluations of propositional logic can be identified with atoms of the Lin\-den\-ba\-um--Tarski algebra of $\models$-equivalent formulas, and any probability mapping is determined by its value on the atoms via the rule (P3). This argument can be used to show that probability functions are exactly the convex combinations of the evaluations. 
A more general proof and precise statement is in Paris \cite{Paris2005}:

\begin{theorem}[see Corollary 4, and Theorem 2 in \cite{Paris2005}]\label{pre:conv}
	A function $B:\gF\to[0,1]$ does not permit a Dutch book if and only if  for every finite $\Gamma\subseteq\gF$, the restriction $B\upharpoonright\Gamma$ is a convex combination of the restrictions of the evaluations $\{ v\upharpoonright\Gamma:\; v\in \Hom(\gF, \TWO)\}$.
\end{theorem}

As formulas are generated by the set of propositional letters $P$, the above restrictions to finite $\Gamma$'s are essentially the same as considering $B$ and the evaluations on the subalgebras $\gF(Q, Cn)$ of $\gF(P, Cn)$ for finite $Q\subseteq P$. We can thus combine Paris' and Ramsey--De Finetti theorems as:

\begin{theorem}\label{thm:pre1}
	Let $(\gF(P, Cn), \models)$ be a classical propositional logic, and $B:\gF\to[0,1]$ a belief function.
	Then the following are equivalent.
		\begin{enumerate}[(A)]\itemsep-2pt
			\item {\bf No Dutch book}:  There is no Dutch book against $B$.
			\item {\bf Probability}: $B$ satisfies the probability axioms (P1)-(P3).
			\item {\bf Convex combination}: For any finite $Q\subseteq P$, $B\restr \gF(Q,Cn)$
			is a convex combination of the functions in $\Hom( \gF(Q,Cn) ,\TWO)$.
			\qed
		\end{enumerate}
\end{theorem}

\paragraph{Two-valued logics with classically working $\land$, $\lor$.}
Results similar to Theorem \ref{thm:pre1} have been carried out in the case of non-classical logics. For two-valued logics the most general result is given by Paris \cite{}, and can be stated as follows. 

\begin{theorem}[Paris \cite{Paris2005}]\label{thm:paris2}
	Let $Cn$ be a set of connectives, which contains the standard $\land$ and $\lor$ (or, these might be derived connectives), 
	and let $\gF = \gF(P,Cn)$ be the formula algebra generated by $P$.
	Let $W$ be a set of $(\land,\lor)$-homomorphisms $w:\gF\to \TWO$, and suppose $\models$ is defined such that
	\begin{equation}
		\phi\models\psi\quad\text{iff}\quad (\forall w\in W)\;
		(w(\phi)=1 \ \ \Rightarrow\ \  w(\psi)=1).
	\end{equation}
	Then the following are equivalent for any $B:\gF\to[0,1]$.
	\begin{enumerate}[(A)]\itemsep-2pt
		\item {\bf No Dutch book}: There is no Dutch book against $B$
		\item {\bf Convex combination}: For any finite $Q\subseteq P$, $B\restr \gF(Q,Cn)$
			is a convex combination of the functions in $\{w\!\upharpoonright \!\gF(Q,Cn):\; w\in W\}$
		\item {\bf Probability}: $B$ satisfies the axioms below.
			\begin{itemize}\itemsep-2pt
				\item[($\mathcal{L}1$)] If $\models\phi$ then $B(\phi)=1$, and
				if $\phi\models$ then $B(\phi)=0$,
				\item[($\mathcal{L}2$)] If $\phi\models\psi$ then $B(\psi)\leq B(\psi)$,
				\item[($\mathcal{L}3$)] $B(\phi\lor\psi)+B(\phi\land\psi)= B(\phi)+B(\psi)$.
				\qed
			\end{itemize}
	\end{enumerate} 
\end{theorem}

\noindent We presented Theorem \ref{thm:paris2} in a more algebraic way than Paris. The condition that each $w\in W$ is a homomorphism into $\TWO$ with respect to $\land$ and $\lor$ amounts to saying that the connectives $\land$ and $\lor$ behave classically. These are precisely the requirements $(T2)$ and $(T3)$ in \cite{Paris2005}. \\

Paris \cite{Paris2005} gives an elementary proof of Theorem \ref{thm:paris2}, and an alternative proof can be given using Choquet's \cite[41.1]{Choquet1954} (see also Bradley \cite{Bradley2016}). Theorem \ref{thm:paris2} applies to a number of well-known propositional logics, for example the standard modal logics $K$, $T$, $S_4$, $S_5$, etc.and to certain paraconsistent logics in which conjunction and disjunction retain their classical interpretation. For similar results concerning two-valued logics we refer to the Dempster--Shafer belief functions, see Jaffray \cite{Jaffray1989}, Shafer \cite{Shafer1976} or Paris \cite{Paris2005}.\\

\paragraph{{\L}ukasiewicz's many-valued logics $\mathL_{k+1}$.} Paris \cite{Paris2005}
proves the analogous result for {\L}ukasiewicz's many-valued logics $\mathL_{k+1}$.
Formulas of $\mathL_{k+1}$ are built up using the standard connectives $\lor$,
$\land$, $\lnot$ and $\to$. Write $\gA = \gA_{k+1}$ for the algebra with the universe 
\begin{equation}
	0,\ \nicefrac{1}{k},\ \nicefrac{2}{k},\ \ldots,\ 1,
\end{equation}
and interpret $\lnot$, $\lor$, $\land$ and $\to$ in $\gA$ by 
the functions $1-a$, $\min\{1, a+b\}$, $\max\{0, a+b-1\}$ and
$\min\{1, 1-a+b\}$ for $a, b\in\gA$. The elements of the universe of $\gA$ are the possible
truth values and the algebraic structure of $\gA$ yields the truth tables of
the logical connectives. Let $W$ be the set of all homomorphisms from $\gF$ into $\gA$.
The consequence $\models$ is defined in the standard way using $W$.

\begin{theorem}[Paris \cite{Paris2005}]\label{thm:luk}
	Let $(\gF(P,Cn), \models)$ be {\L}ukasiewicz's $\mathL_{k+1}$ and
	let $W$ be the set of all homomorphisms from $\gF$ into $\gA_{k+1}$.
	Then the following are equivalent for any $B:\gF\to[0,1]$.
	\begin{enumerate}[(A)]\itemsep-2pt
		\item {\bf No Dutch book}: There is no Dutch book against $B$.
		\item {\bf Convex combination}: For any finite $Q\subseteq P$, $B\restr\gF(Q,Cn)$
			is a convex combination of the functions in $\{w\restr\gF(Q,Cn):\; w\in W\}$.
		\item {\bf Probability}: $B$ satisfies the axioms below.
			\begin{itemize}\itemsep-2pt
				\item[(\L 1)] If $\models\phi$ then $B(\phi)=1$, and
				if $\phi\models$ then $B(\phi)=0$,
				\item[(\L 3)] $B(\phi\lor\psi)+B(\phi\land\psi)= B(\phi)+B(\psi)$.\qed
			\end{itemize}
	\end{enumerate}
\end{theorem}

\noindent There are three other non-classical logics for which the axiomatization of the convex hull is known, and which we recall here.

\paragraph{Kleene's, Priest's and Symmetric logic.} 
Probabilities over Kleene's ``strong logic of indeterminacy'' (KL), Priest's ``logic of paradox'' (LP), and ``symmetric logic'' (SL) are considered in \cite{Williams2016} (for properties of these logics, see \cite{Priest2008}). In each case the set $\gF$ is generated by a non-empty finite set $P$ of propositional variables using the logical connectives $\land$, $\lor$ and $\lnot$.
An \emph{evaluation} (or truth assignment) $h$ assigns \emph{truth values} to propositional variables. Here, we have three the possible truth statuses: $\TRUE$, $\MIDDLEi$, and $\FALSE$. Each evaluation extends to a mapping $h:\gF\to\{\TRUE, \MIDDLEi, \FALSE\}$ 
by the rules given by the \emph{Kleene truth tables} as follows:
\begin{center}
	\begin{tabular}{ c | c  c  c }
	$\land$ & $\TRUE$ & $\MIDDLEi$ & $\FALSE$ \\ \hline
	$\TRUE$ & $\TRUE$ & $\MIDDLEi$ & $\FALSE$  \\ 
	$\MIDDLEi$ & $\MIDDLEi$ & $\MIDDLEi$ & $\FALSE$ \\ 
	$\FALSE$ & $\FALSE$ & $\FALSE$ & $\FALSE$  
	\end{tabular}
	\quad\quad
	\begin{tabular}{ c | c  c  c }
	$\lor$ & $\TRUE$ & $\MIDDLEi$ & $\FALSE$ \\ \hline
	$\TRUE$ & $\TRUE$ & $\TRUE$ & $\TRUE$  \\ 
	$\MIDDLEi$ & $\TRUE$ & $\MIDDLEi$ & $\MIDDLEi$ \\ 
	$\FALSE$ & $\TRUE$ & $\MIDDLEi$ & $\FALSE$  
	\end{tabular}
	\quad\quad
	\begin{tabular}{ c | c  c  c }
	$\lnot$ & $\TRUE$ & $\MIDDLEi$ & $\FALSE$ \\ \hline
	 & $\FALSE$ & $\MIDDLEi$ & $\TRUE$ 
	\end{tabular}
\end{center}
The difference between KL, LP and Symmetric logic is in the definition of their
consequence relations $\models$:
\begin{description}\itemsep-1pt
	\item[{\bf Kleene logic:}] $\phi\models_{KL}\psi$ iff for every evaluation $h$ we have
	\begin{align}
		\text{if } h(\phi)=\TRUE,   &\text{ then } h(\psi) = \TRUE.
	\end{align}
	\item[{\bf LP:}] $\phi\models_{LP}\psi$ iff for every evaluation $h$ we have
	\begin{align}
		\text{if } h(\phi)=\TRUE\text{ or }\MIDDLEi,  &\text{ then } 
		h(\psi) = \TRUE\text{ or }\MIDDLEi.
	\end{align}
	\item[{\bf Symmetric logic:}] $\phi\models_{SL}\psi$ iff for every evaluation $h$ 
	we have	
	\begin{align}
		\text{if } h(\phi)=\TRUE,   &\text{ then } h(\psi) = \TRUE; \text{ and} \\
		\text{if } h(\phi)=\MIDDLEi, &\text{ then } h(\psi) = \TRUE\text{ or }\MIDDLEi.
	\end{align}
\end{description}
In KL the excluded middle $\phi\lor\lnot\phi$ is not a tautology (in fact, this logic has no tautologies at all). LP is a paraconsistent logic, where $\phi\land\lnot\phi$ is not contradictory. SL has both features. The characterization of probabilities in symmetric logic and that of KL and LP with the truth values described here were given in \cite{GSGyW2022a} and is as follows. 

\begin{theorem}[\cite{GSGyW2022a}]\label{thm:sym}
	Let $(\gF, \models)$ be a symmetric logic and let $W$ be the set of all Kleene evaluations. The following are equivalent for any $B:\gF\to[0,1]$.
	\begin{enumerate}[(A)]\itemsep-2pt
		\item {\bf No Dutch book}: There is no Dutch book against $B$.
		\item {\bf Convex combination}: For any finite $Q\subseteq P$, $B\restr\gF(Q,Cn)$
			is a convex combination of the functions in $\{w\restr\gF(Q,Cn):\; w\in W\}$.
		\item {\bf Probability}: $B$ satisfies the axioms below.
			\begin{itemize}\itemsep-2pt
				\item[(SL1)] If $\phi\models\psi$ then $B(\phi) \leq B(\psi)$,
				\item[(SL2)] $B(\lnot \phi) = 1-B(\phi)$,
				\item[(SL3)] $B(\phi\lor \psi) = B(\phi) + B(\psi) - B(\phi\land \psi)$,
				\item[(SL4)] $B(\phi)=B(\phi\land\psi) + B(\phi\land\lnot\psi) - 
					B(\phi\land \lnot\phi\land \psi\land \lnot\psi)$. \qed
			\end{itemize}
	\end{enumerate}
\end{theorem}

\noindent In case of KL and LP the similar theorem holds with the only modification that (SL1) should be replaced with
\begin{quote}
	(KLP1)\ If $\phi\models\psi$ and $\lnot\psi\models\lnot\phi$, then $B(\phi) \leq B(\psi)$,
\end{quote}
\noindent The reason is that the semantic consequence relation of KL, LP, and SL are interdefinable:
\begin{enumerate}[(i)]\itemsep-2pt
	\item $\phi\models_{SL}\psi$ if and only if ($\phi\models_{KL}\psi$ and $\lnot\psi\models_{KL}\lnot\phi$).
	\item $\phi\models_{LP}\psi$ if and only if ($\lnot\psi\models_{KL}\lnot\phi$).
	\item $\phi\models_{SL}\psi$ if and only if ($\phi\models_{LP}\psi$ and $\lnot\psi\models_{LP}\lnot\phi$).
\end{enumerate}	
\medskip

\section{General patterns} \label{sec:genpat}

In the previous theorems we had a hidden assumption that the truth values which define the logical consequence relation are actual real numbers. Williams introduces the term ``cognitive load'': a cognitive load of a truth value is the supposed ``ideal cognitive state'' associated with it; in other words, it is the degree of belief an agent should invest in a proposition having that truth value. In the classical case, cognitive loads directly correspond to truth values $1$ (true) and $0$ (false), while in the general case the cognitive load function $e$ is an arbitrary function from the truth values into $[0,1]$. For any logical evaluation $w$ we can speak of ``its'' \textit{cognitive evaluation} $e\circ w: \gF\to [0,1]$. Williams's idea is, in the context of some logic, to inquire about the convex combinations of something else than valuations. The reason for this is that two different logics, defined on the same language and having the same set of truth values, may give rise to exactly the same set of valuations. And yet, for example due to how the consequence relation differs between the two logics, the epistemic status of these valuations might be different. According to both Williams \cite{Williams2016} and Bradley \cite{Bradley2016} the logics are ``cognitively loaded'', in that each truth value has ``its'' cognitive load \citep[255]{Williams2016}. For the purposes of stating the formal results, in the upcoming sections we treat the cognitive load function to be a definitional element of the given logic.

Formally, let $\gA$ be any algebra (of truth values). In the examples above: in the two-valued case $\gA$ the lattice (Boolean algebra) $\TWO$, in the case of {\L}ukasiewicz's $\mathL_{k+1}$ it was the $k+1$-element algebra $\gA_k$, and in the case of Kleene, Priest and Symmetric logic, it was a $3$-element algebra, the operations of which were defined by the Kleene truth tables. The elements of $\gA$ are abstract truth values. Each such truth value has a numeric value given by the cognitive load function $e:\gA\to [0,1]$. For a given set $P$ of propositional letters, the formula algebra $\gF(P)$ is freely generated by $P$. Evaluations of a logic whose algebra of truth values is $\gA$ is the set of homomorphisms $\Hom(\gF, \gA)$ between the formula algebra and the algebra of truth values. For each such homomorphism $w$, the function (cognitive evaluation) $e\circ w$ maps formulas into $[0,1]$
\begin{align}
		\gF \overset{w}{\longrightarrow} \gA\overset{e}{\longrightarrow}[0,1]\,.
\end{align}
To axiomatize probabilities (not Dutch bookable functions) over a given logic we had to axiomatize the convex hull of the functions $e\circ w$ for $w\in \Hom(\gF, \gA)$. Apart from the evaluations the other main component of logics was the consequence relation $\models$. The $3$-valued logics above show that it may happen that different logics have the same set of evaluations. Observe that $\models$ has absolutely no role in what is the convex hull of the cognitive evaluations. Let us take a look at the shape of the axioms in the previous section. We can make a distinction between two types of axioms. 

The first type of axioms connects the logical consequence relation $\models$ and the belief function $B$, e.g. (P1), (P2), (SL1). Let us call an axiom of this type a logical axiom. We will see later on that the main purpose of logical axioms is to make sure that equivalent formulas get the same probability. More precisely, call two formulas $\varphi, \psi\in\gF$ \emph{equivalent}, $\varphi\sim\psi$, if for any evaluation $w\in \Hom(\gF, \gA)$ we have $w(\varphi)=w(\psi)$. Now take any function $B:\gF\to[0,1]$. If $B$ belongs to the convex hull of cognitive evaluations, then $B(\varphi)=B(\psi)$ must hold for equivalent $\varphi$ and $\psi$. In other words, the property 
\begin{align}
	\varphi\sim\psi\quad\Longrightarrow\quad B(\varphi)=B(\psi)\label{eq:nec}
\end{align}
is \emph{necessary} for $B$ to belong to the convex hull of cognitive evaluations. In the examples above (and, in fact, most of the logics, e.g. Blok--Pigozzi algebraizable logics) equivalence can be expressed by $\models$. for example, in classical propositional logic $\varphi$ and $\psi$ are equivalent, if $\varphi\models\psi$ and $\psi\models\varphi$. Now, axioms such as (P1) and (P2) basically gives: (i) if $\varphi$ is equivalent to ``true'', then $B(\varphi) = $ ``$B(true)$'' $=1$; (ii) if $\varphi$ and $\psi$ are equivalent, then $B(\varphi)=B(\psi)$. In the case of KL or LP expressing equivalence of $\varphi$ and $\psi$ is more complicated: $\varphi$ and $\psi$ are equivalent if and only if
\begin{align}
	\varphi\models\psi,\quad\text{ and }\quad \psi\models\varphi,\quad\text{ and }\quad \lnot\varphi\models\lnot\psi,\quad\text{ and }\quad \lnot\psi\models\lnot\varphi\,.
\end{align}
Hence the unusual axiom (KLP1) above. Since \eqref{eq:nec} is a necessary condition, to avoid further definitions and discussion on how exactly equivalence is expressible in terms of $\models$, instead of these so called logical axioms we will simply assume that \eqref{eq:nec} is an axiom, and we will not make use of any explicit mention of $\models$. Note that in the usual cases the definition of $\models$ depends on what the set $P$ is, and thus an axiom like (P1) or (P2) in fact depend on $P$ (though, the ``form'' does not change for different $P$'s).  

The second type of axioms, e.g. (P3), (SL2)-(SL4), stipulate equalities (in general, inequalities) between certain expressions, and there is no mention of $\models$. Call this type of axioms non-logical. Take, for example, $B(\varphi\lor\psi) = B(\varphi)+B(\psi)-B(\varphi\land\psi)$. We could think of such an axiom as if $B$ was defined not only on the formulas, but also on ``formal expressions'' of the form $\varphi\lor\psi = \varphi + \psi - (\varphi\land\psi)$. Clearly the symbols $+$, $-$ are not logical connectives, but we might consider such expressions as some sort of extensions of the formulas. We will make this idea precise later on. The main point we would like to mention is that for any given logics above, the set of axioms are the same for any choice of $P$. Thus, with a fixed algebra of truth values, if the logics corresponding to different $P$'s  are treated as being instances\footnote{For example, we can think of  ``classical propositional logic'' as the collective word for any particular classical propositional logic generated by some set of propositional letters $P$. This idea is made precise in \cite{LogicFamilies}.} of the same ``general'' logic, then it is only the ``general'' logic on which the axioms of probability depend.

\paragraph{Lindenbaum--Tarski algebras.}
Fix the algebra $\gA$ of truth values, and the set $P$ of propositional letters that generate the formula algebra $\gF(P)$. The equivalence relation $\varphi\sim\psi$ defined by 
\begin{align}
	\varphi\sim\psi\quad\Longleftrightarrow\quad \big(\forall w\in\Hom(\gF, \gA)\big)\ w(\varphi)=w(\psi)
\end{align}
is a congruence relation\footnote{For each homomorphism $w$ the kernel $\ker(w) = \{ (\phi,\psi):\; w(\phi)=w(\psi)\}$ is a congruence; $\sim$ equals $\bigcap_{w}\ker(w)$; and the intersection of congruences is a congruence.} of $\gF$. Therefore the quotient structure $\gF/_{\sim}$ exists. This quotient structure is the Lindenbaum--Tarski algebra of the logic generated by $P$ and we use the notation $\gL(P) = \gF(P)/_{\sim}$ (we may omit the mention of $P$ if it is clear from the context). Let us remark that $\gL$ is a free algebra in the variety generated by $\gA$. Also, $\gL$ can be embedded into $\gA^{|\Hom(\gF,\gA)|}$ by the mapping $\phi/_{\sim} \mapsto \< w(\phi): w\in \Hom(\gF,\gA)\>$. As each homomorphism $w:\gF\to\gA$ is 
determined by its values $w\restr P$ on the generator elements, there are exactly $|\gA|^{|P|}$ homomorphisms, hence $\gL$ is finite: $|\gL|\leq |\gA|^{|\gA|^{|P|}}$.

\section{The general characterization theorems}\label{sec:new}

We start with setting the stage. Throughout this section $Cn$ is a \emph{finite} set of logical connectives (algebraic operations), and $\gA$ is a \emph{finite} $Cn$-type algebra, the algebra of truth values (the universe of this algebra will be denoted by $A$). $\gF(P)$ is the formula algebra freely generated by $P$ in similarity type $Cn$. We write $\sim_P$ for the congruence
\begin{align}
	\varphi\sim_{P}\psi\quad\Longleftrightarrow\quad \big(\forall w\in\Hom(\gF(P), \gA)\big)\ w(\varphi)=w(\psi),
\end{align}
and write $\gL(P)$ for $\gF(P)/_{\sim}$. Whenever it does not cause confusion we may omit the parameter $P$ from $\gF(P)$, $\sim_{P}$, etc. The cognitive load is an arbitrary function $e:\gA\to [0,1]$. We write
\begin{align}
	\Hom_e(\gF,\gA) = \big\{ e\circ w:\; w\in \Hom(\gF,\gA) \big\}\,
\end{align}
for the set of cognitive evaluations.

\subsection{Birkhoff's equational logic extended}\label{subsec:birkhoff}

Equational logic (see e.g. \cite[II.14]{BurrisSankappanavar}) gives a precise account on deriving identities from other identities. We do not recall here the definitions, but an illustrative example is the case of Boolean algebras: every identity that holds in all Boolean algebras, e.g. $(x\lor y)\lor z = x\lor (y\lor z)$, can be derived from the Boolean axioms. These axioms are formulated by using two variables only. We extend this calculus to probability axioms: from certain axioms we can infer to other axioms. An obvious, not so interesting example could be: from the axiom $B(\lnot\varphi) = 1-B(\varphi)$ we can infer to $B(\varphi) = 1-B(\lnot\varphi)$; or to $B(\varphi\land\psi) = 1-B(\lnot\varphi\lor\lnot\psi)$ 
(say, in propositional logic). 
In what follows, we make this idea precise. 

A remark about axioms is in order here. By an axiom, say, $B(\varphi\lor\psi) \leq B(\varphi)+B(\psi)$, we understand, that this inequality holds for any choice of $\varphi$ and $\psi$. As customary in defining derivation systems, we pick a set $X$ of \emph{formula variables}, and phrase the axiom as 
$B(x\lor y) \leq B(x)+B(y)$. One can substitute concrete formulas in place of formula variables, to get an instance of the axiom.

\begin{definition}
	For a set $X$ of (formula) variables, the set of \emph{formal expressions} over $X$ is defined 
	as the smallest set (with respect to inclusion) $\gT(X)$ such that
	\begin{itemize}\itemsep-2pt
		\item every $\varphi\in\gF(X)$,
		\item if $t_1$, and $t_2$ are in $\gT(X)$, then so are $t_1+t_2$ and $r\cdot t_1$
		for $r\in \mathbb{R}$.
	\end{itemize}
	For a formal expression $t\in \gT(X)$ and $r\in\mathbb{R}$ we call $t \geq r$ a \emph{formal inequality}.
	The set of formal inequalities is denoted by $\gT^{\geq}(X)$.
\end{definition}

\noindent For a finite set $I$ of indices, and $\varphi_i\in \gF(X)$, 
$\alpha_i, r\in \mathbb{R}$ ($i\in I$), 
$\sum_{i\in I}\alpha_i\varphi_i$ is a formal expression, 
and $\sum_{i\in I}\alpha_i\varphi_i\geq r$ is a formal inequality.
It is clear that every formal expression (inequality) is of this form. 
Next, we extend functions $B:\gF(P)\to[0,1]$ to formal expressions.

\begin{definition}\label{def:extf}
	For a finite set $I$ of indices, and $\varphi_i\in \gF(P)$, $\alpha_i\in \mathbb{R}$ ($i\in I$), 
	and $B:\gF(P)\to[0,1]$ we define
	\[
		B\big( \sum_{i\in I}\alpha_i\varphi_i \big) = \sum_{i\in I}\alpha_iB(\varphi_i)\,.
	\]
	We say that a function $B:\gF(P)\to [0,1]$ \emph{validates} a formal inequality $t\geq r\in\gT^{\geq}(X)$, 
	if \emph{for all} $\varphi_1$, $\ldots$, $\varphi_n\in \gF(P)$ we have
	\[
		B\big( t(\varphi_1, \ldots, \varphi_n)  \big) \geq r\,.
	\]
	Here $t(\varphi_1, \ldots, \varphi_n)$ means that we substituted formulas 
	$\varphi_i\in\gF(P)$ in place of the variables $x_i\in X$. 
	We write $B\models t\geq r$ to denote that $B$ validates $t \geq r$. For $\Sigma\subseteq\gT^{\geq}(X)$, 
	we write $B\models \Sigma$ to denote that $B\models t\geq r$ for every $t\geq r\in\Sigma$.
\end{definition}
	
\noindent Generally, the non-logical axioms of probability are of the form 
$B\models t\geq r$ for some $t\in \gT(X)$, and $r\in\mathbb{R}$.
We use the standard conventions that $t\leq r$ means $-t\geq -r$, 
$t_1\leq t_2$ means $t_1-t_2\leq 0$, and $t_1=t_2$ denotes the two inequalities $t_1\leq t_2$ and $t_2\leq t_1$. For example, the axiom (P3) can be written as
\[
	B\models (x\lor y) = x+y-(x\land y)\,,
\] 
and this means that whatever formulas (from the domain of $B$) we substitute in place of the formula 
variables $x$ and $y$, after the substitution and applying $B$ we get an equality that holds. 
It is important to note that we do not allow strict inequalities in formal inequalities.

\begin{definition}
	For $B:\gF(P)\to[0,1]$ we define 
	\[
		\Id_{B}^{\geq}(X) = \big\{ t\geq r\in\gT^{\geq}(X):\; B\models t\geq r \big\}\,.
	\]
	Similarly, 
	\[
		\Id^{\geq}(X) = \big\{ t\geq r\in\gT^{\geq}(X):\; (\forall P)(\forall B:\gF(P)\to[0,1])\ 
		 B\models t\geq r  \big\}\,.
	\]
	For every $B$, $B\models\Id_{B}^{\geq}(X)$ and $B\models\Id^{\geq}(X)$.
\end{definition}

\noindent Formal inequalities can be added, subtracted, multiplied in the obvious way. 
For example, $t_1\geq r_1$ added to $t_2\geq r_2$ is $t_1+t_2\geq r_1+r_2$. 
Observe that for $t(x), p \in\gT(X)$ the substituted expression $t(x/p)$ also belongs to $\gT(X)$. 
Further, if $B\models t\geq r$, then for any $p\in \gT(X)$ we have $B\models t(x/p)\geq r$. 
Similar remark goes for replacing variables with other variables. The operation of substitution and replacement of variables correspond to applying endomorphisms $\alpha\in\Hom(\gT(X), \gT(X))$ to the terms $t\in\gT(X)$.

\begin{observation}\label{obs:id}
	The sets $\Id_{B}^{\geq}(X)$ and $\Id^{\geq}(X)$ contain
	\begin{itemize}\itemsep-2pt
		\item $t=t$ for every $t\in\gT(X)$,
		\item $\varphi\geq 0$ and $\varphi\leq 1$ for $\varphi\in\gF(X)$,
	\end{itemize}
	and are closed under the following operations:
	\begin{itemize}\itemsep-2pt
		\item Addition, subtraction and scalar multiplication by non-negative real numbers.
		\item Endomorphisms $\alpha\in \Hom(\gT(X), \gT(X))$: if $t\geq r$ belongs to
		the set, then so does $\alpha(t)\geq r$. 
		\item If $t_1\leq t_2$, and $t_2\leq t_3$ belongs to the set, then so does
		$t_1\leq t_3$. 
		\item Applying commutativity, associativity and distributivity of $+$, $-$, and $\cdot$, for example,
		if $t_1+t_2\geq r$ is in the set, then so is $t_2+t_1\geq r$.
		\item If $t\geq r$ belongs to the set, then so is $t\geq r'$ for any $r'\leq r$. 
	\end{itemize}
\end{observation}

\noindent For a set $\Sigma\subseteq\gT^{\geq}(X)$, we write $B\models \Sigma$ if $B\models t\geq r$ 
for each $t\geq r\in \Sigma$. Validity of formal inequalities induce a semantic consequence 
relation between sets of inequalities and inequalities as follows.	

\begin{definition}
	Let $\Sigma\cup\{t\geq r\}\subseteq\gT^{\geq}(X)$. We say that $\Sigma\models t\geq r$
	if for every $B:\gF(P)\to[0,1]$ whenever $B\models \Sigma$, then $B\models t\leq r$ as well.
\end{definition}

\begin{theorem}[Compactness of validity]
	For any $\Sigma\subseteq\gT^{\geq}(X)$ the following are equivalent
	\begin{enumerate}[(A)]\itemsep-1pt
		\item There is $B:\gF(P)\to[0,1]$ such that $B\models\Sigma$.
		\item For each finite $\Gamma\subseteq\Sigma$ there is $B_{\Gamma}:\gF(P_{\Gamma})\to[0,1]$
		such that $B_{\Gamma}\models \Gamma$.
	\end{enumerate}
\end{theorem}	
\begin{proof}
	(A)$\Rightarrow$(B) is straightforward. The direction (B)$\Rightarrow$(A) is an ultralimit argument.
	Assume (B). Without loss of generality, we can suppose 
	that the sets $P_{\Gamma}$ are the same for every $\Gamma$, because
	satisfaction of an inequality is invariant under permutations of variables, and each formula can use 
	finitely many variables only, thus a countable infinite set $P$ of propositional variables 
	is always enough. 
	Let $I$ be the set of finite subsets of $\Sigma$, and take a non-trivial ultrafilter $U$ over $I$. 
	Define $B:\gF(P)\to [0,1]$ as the ultralimit $B(\varphi) = \lim_{U}B_{i}(\varphi)$.
	As $[0,1]$ is compact and Hausdorff, this ultralimit always exists, and thus for each $\varphi$
	there is $B(\varphi)\in\mathbb{R}$ such that for every $\varepsilon>0$,
	$\big\{ i\in I:\; |B_i(\varphi) - B(\varphi)|<\varepsilon \big\}\in U$. This $B$ extends
	to formal expressions as in Definition \ref{def:extf}. Using linearity of the 
	ultralimit, for any formal expression $t$, $B(t) = \lim_{U}B_i(t)$. For
	a formal inequality $t\geq r\in \Sigma$, $B\models t\geq r$ holds whenever 
	$\{i:\; B_i\models t\geq r\}\in U$. To achieve this, it is enough to chose the ultrafilter such 
	that $\{i:\; t\geq r\in i\}\in U$. This is possible, as the system of these latter sets have the 
	finite intersection property, and hence extend to an ultrafilter.
\end{proof}

\noindent Compactness of semantic consequence does not hold, in general. For, take
$\Sigma = \{\varphi\geq 1-\nicefrac{1}{n}:\; n\in\mathbb{N} \}$. Then $\Sigma\models \varphi=1$, but 
for no finite subset $\Gamma\subseteq\Sigma$ it holds that $\Gamma\models \varphi=1$.

\begin{definition}\label{def:DSigma}
	For $\Sigma\subseteq\gT^{\geq}(X)$ we define the \emph{deductive closure} $D(\Sigma)$ of $\Sigma$ as 
	the smallest subset of $\gT^{\geq}(X)$ such that 
	\begin{itemize}\itemsep-1pt
		\item $\varphi\geq 0$, $\varphi\leq 1$, and	$t=t\in D(\Sigma)$ for every $\varphi\in\gF(X)$, 
		$t\in \gT(X)$.
		\item $D(\Sigma)$ is closed under addition, subtraction and scalar multiplication 
		by non-negative real numbers.
		\item $D(\Sigma)$ is closed under endomorphisms $\alpha\in\Hom(\gF(X),\gF(X))$
		(that is, closed under substitutions and replacements of variables).
		\item If $t_1\leq t_2$, and $t_2\leq t_3$ are in $D(\Sigma)$, then $t_1\leq t_3\in D(\Sigma)$.
		\item $D(\Sigma)$ is closed under commutativity, associativity and distributivity of $+$, 
		$-$, and $\cdot$. 		
		\item If $t\leq r\in D(\Sigma)$, then $t\leq r'\in D(\Sigma)$ for any $r\leq r'$. 
		\item If $\sum\alpha_it_i\leq r$ is in $D(\Sigma)$ for $\alpha_i\geq 0$ but $r<0$, then 
		$D(\Sigma)=\gT^{\geq}(X)$.
	\end{itemize}
	The last item expresses that if there is one contradictory inequality in the system, then every inequality 
	is in the system (principle of explosion).
\end{definition}

\begin{proposition}\label{prop:d}
	For any $\Sigma\subseteq \gT^{\geq}(X)$, if $B\models \Sigma$, then $B\models D(\Sigma)$.
\end{proposition}
\begin{proof}
	If $B\models\Sigma$, then $\Sigma\subseteq \Id_{B}^{\geq}(X)$. But $\Id_{B}^{\geq}(X)$ is
	closed under all the closure properties in the definition of $D(\Sigma)$ (cf. Observation \ref{obs:id}).
	Therefore, $D(\Sigma)\subseteq \Id_{B}^{\geq}(X)$, and thus $B\models D(\Sigma)$.	
\end{proof}

\begin{corollary}\label{cor:sound}
	For $\Sigma\cup\{t\geq r\}\subseteq \gT^{\geq}(X)$, if $t\geq r\in D(\Sigma)$, 
	then $\Sigma\models t\geq r$.
\end{corollary}
\begin{proof}
	Take any $B\models \Sigma$ (if there is no such, then $\Sigma\models t\geq r$ by definition).
	By Proposition \ref{prop:d}, $B\models D(\Sigma)$, and thus $B\models t\geq r$. As $B$ was
	arbitrary, we obtain $\Sigma\models t\geq r$.
\end{proof}

The deductive closure can be turned into a formal inference system as follows. 

\begin{definition}
	Let $\Sigma\cup\{t\geq r\}\in \gT^{\geq}(X)$. 
	We define $\Sigma\proves t\geq r$ (reads $\Sigma$ proves/derives $t\geq r$) 
	if there is a finite sequence of inequalities $p_1\geq r_1$, $\ldots$, $p_n\geq r_n$
	from $\gT^{\geq}(X)$ such that each $p_i\geq r_i$ either
	\begin{itemize}\itemsep-1pt
		\item belongs to $\Sigma$, or
		\item of the form $\varphi\geq 0$, $\varphi\leq 1$ or $t=t$ for $\varphi\in\gF(X)$, $t\in\gT(X)$, or
		\item is a result of applying any of the closure rules from Definition \ref{def:DSigma} to 
		previous inequalities in the sequence,
	\end{itemize}
	and the last inequality is $t\geq r$. Such a sequence is called a \emph{formal deduction} of $t\geq r$
	from $\Sigma$
\end{definition}

\begin{definition}
	$\Sigma\subseteq\gT^{\geq}(X)$ is \emph{semantically consistent} if there is $B\models\Sigma$, 
	and semantically inconsistent otherwise. We say that $\Sigma$ is \emph{syntactically consistent}
	if $D(\Sigma)\neq\gT^{\geq}(X)$, or equivalently, $\Sigma$ does not prove any contradictory 
	inequality (cf. the last item in Def. \ref{def:DSigma}); otherwise it is syntactically inconsistent.
\end{definition}

It is straightforward that a syntactically inconsistent $\Sigma$ is also semantically inconsistent. Equivalently, a semantically consistent $\Sigma$ is syntactically consistent.

\begin{theorem}[Soundness]
	Given $\Sigma\cup\{t\geq r\}\in \gT^{\geq}(X)$, 
	if $\Sigma\proves t\geq r$, then $\Sigma\models t\geq r$.
\end{theorem}
\begin{proof}
	Immediate from Corollary \ref{cor:sound} because $\Sigma\proves t\geq r$ iff $t\geq r\in D(\Sigma)$.	
\end{proof}

So far the underlying propositional logic did not play any role. Consider the variety $\HSP(\gA)$. 
By Birkhoff's theorem there is a set $\Gamma$ of identities of the form $t_1=t_2$ for $t_1, t_2\in\gF(X)$ 
for a suitable set $X$ of variables (and thus they are formal inequalities) such that members of the variety are exactly the algebras satisfying all identities in $\Gamma$. If $w\in\Hom(\gF(P), \gA)$ and 
$t_1=t_2\in \Gamma$, then $w(t_1)=w(t_2)$ for any substitution of formulas in place of the variables in $X$, that is, $w\models t_1=t_2$.
According to the general scheme of things in this paper, we are interested in axiomatizing functions $B$
that lie in the convex hull of the cognitive evaluations. Thus, the expressions $t_1=t_2$ in $\Gamma$ can be considered as potential axioms.
Take now any $\varphi, \psi\in\gF(X)$ that are equivalent $\varphi\sim_{X}\psi$. Then the identity $\varphi=\psi$ is true in the variety, and by Birkhoff's calculus (cf. \cite[Theorem 14.19]{BurrisSankappanavar}), $\varphi=\psi$ is derivable from $\Gamma$. Thus, if $B\models\Gamma$, 
then $B(\varphi)=B(\psi)$. 
It follows that the ``logical axiom'' (if $\varphi\sim\psi$ then $B(\varphi)=B(\psi)$) can be replaced by the ``non-logical axioms'' ($B(t_1)=B(t_2)$ for $t_1=t_2\in\Gamma$). If the variety
is finitely based, that is, $\Gamma$ can be taken finite (such is the case in classical propositional logic, or Kleene's, Priest's, and Symmetric logic, etc.), then this means adding finitely many non-logical axioms in place of the single logical axiom in the characterization of the convex hull of the cognitive evaluations.

Let $X$ be a countable set of variables. 
We define the formal inequalities imposed by the underlying logic as the set
\[
	\Id_{\cL} = \big\{ \varphi=\psi:\; \varphi,\psi\in\gF(X), \text{ and } \gA\models \varphi=\psi \big\}.
\]

Recall that a syntactically inconsistent $\Sigma$ is semantically inconsistent; equivalently, a semantically consistent $\Sigma$ is syntactically consistent. The next theorem states the converse of this implication under some additional assumptions. 

In the proof we make use of Farkas' lemmas (see \cite{Dantzig}). For a real matrix $A\in\mathbb{R}^{m\times n}$ and vector $\vec r\in\mathbb{R}^{m}$ exactly one of the following alternatives holds. 
Either there exists $\vec x\in\mathbb{R}$ such that $A\cdot\vec x\leq \vec r$ and $\vec x\geq 0$; or there exists $\vec y\in\mathbb{R}^m$ such that $\vec y\geq 0$, $\vec y\cdot A\geq 0$, and 
$\vec y\cdot \vec r<0$.  We make use of a strict version (Ville's alternative, \cite[p.139]{Dantzig}, \cite[Thm.2.10]{Gale}) as well: Either there exists $\vec x\in\mathbb{R}$ such that $A\cdot\vec x < \vec r$ and $\vec x\geq 0$; or 
there exists $\vec y\in\mathbb{R}^m$ such that $\vec y\cdot\vec r\leq 0$, 
and $\vec y\cdot A\geq 0$.

\begin{theorem}[Finite completeness]
	Let $X$ be a countable, and $\Sigma\subseteq\gT^{\geq}(X)$ a finite set. 
	\begin{enumerate}[(a)]\itemsep-1pt
	\item $\Sigma\cup\Id_{\cL}$ is syntactically consistent if and only if it is semantically consistent.
	\item If $\Sigma\cup\Id_{\cL}\models t\leq c$, then $\Sigma\cup\Id_{\cL}\proves t\leq c$.
	\end{enumerate}
\end{theorem}
\begin{proof}
	(a)
	Let $X_0$ be the finite subset of $X$ that consists of the variables occurring in formulas in $\Sigma$. 
	The Lindenbaum--Tarski algebra $\cL(X_0)$ is finite, because $\gA$ is finite. 
	Substitute in all possible
	ways elements of $\cL(X_0)$ into the variables of the formal inequalities in $\Sigma$, 
	and denote this finite set by $\Sigma_0$. Up to $\sim_{X_0}$-equivalence, every formal inequality in
	$\Sigma_0$ can be written as $\sum_{i}\alpha_i\varphi_i\leq r$, where $\varphi_\in\gL(X_0)$. 
	Consider the matrix which has $|\gL(X_0)|$ columns and for each formal inequality 
	$\sum_{i}\alpha_i\varphi_i\leq r$ in $\Sigma_0$ it contains a row consisting of the $\alpha_i$'s. Let 
	$\vec{r}$ be the vector containing the $r$'s of these inequalities at the corresponding coordinates.
	Extend the matrix with the identity matrix, and extend $\vec{r}$ by corresponding $1$'s. Denote 
	this matrix by $A$. 
	
	If there is $\vec{x}$ with all non-negative coordinates ($\vec{x}\geq 0$)
	solving the system of inequalities $A\vec{x}\leq \vec{r}$, then 
	we can define $B:\gF(X_0)\to[0,1]$ as follows. For $\psi\in\gF(X_0)$ let $\varphi_i$ be the
	element of $\gL(X_0)$ which is the representant of the class $\psi/_{\!\sim_{X_0}}$. Then put
	$B(\psi) = x_i$. This gives a mapping $B$ for which
	$0\leq B(\varphi)\leq 1$ (by non-negativity of $\vec{x}$, and 
	the identity block of the matrix $A$), and for each
	$t\leq r\in \Sigma$, $B\models t\leq r$, because for any substitution there is a line in the matrix 
	corresponding to exactly that inequality.
	
	By Farkas' lemma, either $A\vec{x}\leq \vec{r}$ has a solution $\vec{x}\geq 0$; or
	there is $\vec{y}\geq 0$ with $\vec{y}A\geq 0$ and $\vec{y}\cdot \vec{r}<0$. This latter option
	says that from $\Sigma$, using substitutions and linear combinations (that is, formal derivations)
	we can derive an inequality $\sum\beta_i\varphi_i\leq c$ such that $\beta_i\geq 0$, but $c<0$.
	But this means that $\Sigma\cup \Id_{\cL}$ is syntactically inconsistent.  
	
	Summing up, either there is $B\models \Sigma\cup\Id_{\cL}$, and thus 
	$\Sigma\cup \Id_{\cL}$ is semantically consistent; or if there is no such $B$,
	and thus $\Sigma\cup \Id_{\cL}$ is semantically inconsistent,  
	then $\Sigma\cup \Id_{\cL}$ is syntactically inconsistent as well. \\
	
	(b) By part (a) we can assume that $\Sigma\cup\Id_{\cL}$ is syntactically (and thus 
	semantically) consistent, as
	otherwise there is nothing to prove. 
	Consider the matrix $A$ from part (a) of the proof, and assume that 
	$\Sigma\cup\Id_{\cL}\not\proves t\leq c$. 
	Consider $t$ as a vector $\vec{t}$ of coefficients of formulas $\varphi\in\cL(X_0)$ just 
	as we constructed the rows of the matrix $A$. 
	
	Suppose $\Sigma\cup\Id_{\cL}\not\proves t\leq c$. Then there is no $\vec{y}$
	such that 
	\begin{align}
		\vec{y}A=\vec{t}\quad \text{ and }\quad \vec{y}\cdot\vec{r}\leq c. \label{eq:yAtc}
	\end{align}
	Extend the matrix $A$ with an extra row $-\vec{t}$, and 
	extend the vector $\vec{r}$ with an extra coordinate $-c$. Let us
	denote these by $(A|-\vec{t})$ and $(\vec{r}|-c)$. Condition \eqref{eq:yAtc}
	is equivalent to that there is no $\vec{z}$ such that the last coordinate of $\vec{z}$
	is not $0$, and $\vec{z}(A|-\vec{t}) = 0$, while $\vec{z}(\vec{r}|-c)\leq 0$. (This is an 
	equivalent condition, because if there was
	such a $\vec{z}$ but with $0$ as the last coordinate, then $\Sigma\cup\Id_{\cL}$ would be 
	inconsistent). By Farkas lemma (the version for strict inequalities), 
	there exists then $\vec{x}\geq 0$
	such that $(A|-\vec{t})\vec{x} < (\vec{r}|-c)$. It follows that 
	$\vec{t}\cdot\vec{x} > c$ for this $\vec{x}$. Now, 
	$\vec{x}$ gives rise to a function $B:\gF(X_0)\to[0,1]$
	(cf. the end of the proof of part (a)) such that $B\models\Sigma\cup\Id_{\cL}$ 
	but $B\not\models t\leq c$.
\end{proof}

Completeness of $\models$ with respect to $\proves$ does not hold in general.
Let $x\in X$ and take $\Sigma = \{ x\geq 1-\nicefrac{1}{n}:\; n\in\mathbb{N}\}$. Then 
$\Sigma\models x=1$, but $\Sigma\not\proves x=1$ can be shown by induction on the complexity 
of finite formal derivations.

\begin{definition}
	We define the set $\CH(X)\subseteq \gT^{\geq}(X)$ of formal inequalities that are 
	satisfied by all the cognitive evaluations $\Hom_e(\gF(P),\gA)$:
	\[
		\CH(X) = \{ t\geq r\in\gT^{\geq}(X):\; \forall P\ \forall v\in\Hom_e(\gF(P), \gA)\ \ 
		v\models t\geq r\}\,.
	\]
	Intuitively, $\CH(X)$ is the set of all possible axioms in variables $X$ 
	that could characterize the 
	convex hull of $\Hom_e(\gF(P),\gA)$: since every cognitive evaluation 
	$v\in \Hom_e(\gF(P), \gA)$ validates every member of $\CH(X)$, the same is true 
	for all convex combinations of such cognitive evaluations. 
\end{definition}

\noindent If an inequality $t\geq r\in \gT^{\geq}(X)$ is validated by all cognitive homomorphisms in
$\Hom_e(\gF(X), \gA)$, then the same is true for cognitive homomorphisms in $\Hom_{e}(\gF(P),\gA)$
for any set $P$. Therefore, we could have defined $\CH(X)$ as
\[
	\CH(X) = \{ t\geq r\in\gT^{\geq}(X):\; \forall v\in\Hom_e(\gF(X), \gA)\ \ 
		v\models t\geq r\}\,.
\]

\subsection{Separable homomorphisms}\label{sub:sep}

One of the useful tools for the axiomatization of the convex hull of the evaluations available in some logics is the ``separation of evaluations''.
In propositional logic, evaluations can be identified with atoms of the Lin\-den\-ba\-um--Tarski algebra of $\models$-equivalent formulas, and any probability mapping is determined by its value on the atoms via the rule (P3). This argument can be used to show that probability functions are exactly the convex combinations of the evaluations. If, for example, $P=\{p,q\}$, then the atoms of the Lindenbaum--Tarski algebra are $p\land q$, $p\land\lnot q$, $\lnot p\land q$, and $\lnot\land\lnot q$. Any mapping given on these atoms extend additively to the entire Boolean algebra. The atoms and the evaluations are in a one-to-one correspondence: by writing $p^1=p$ and $p^0=\lnot p$ we have
\[
	v\big(p^{w(p)}\land q^{w(q)}\big) = 
	\begin{cases}
		1 & \text{ if } v=w,\\
		0 & \text{ otherwise.}
	\end{cases}
\]
Something similar happens in the case of {\L}ukasiewicz's $k$-valued logics (see Paris \cite{Paris2005}). There, for any truth value $i$ there is a formula $\psi_i(p)$ such that
\[
	w\big( \psi_i(p) \big) = \begin{cases}
		1& \text{ if } w(p) = \nicefrac{i}{k},\\
		0 & \text{ otherwise.}
	\end{cases}
\]
Then letting $\chi_w = \psi_{i_1}(p_1)\land\ldots\land \psi_{i_n}(p_n)$ for $\nicefrac{i_r}{k}=w(p_r)$, we get
\[
	v\big( \chi_w \big) = \begin{cases}
		1& \text{ if } v=w, \\
		0 & \text{ otherwise.}
	\end{cases}
\]
The next theorem shows that this feature of a logic is always sufficient for an axiomatization of the convex hull of the evaluations. 

\begin{theorem}\label{thm:m1}
	Let $P$ be a finite set. Suppose that for each $w\in\Hom_e(\gF(P), \gA)$
	there are $\alpha_i^w\in \mathbb{R}$, and $\varphi_i^w\in\gF(P)$ 
	for $i\in I_w$ such that for every $v\in \Hom_e(\gF(P), \gA)$ we have
	\[
		\sum_{i\in I_w}\alpha_i^w v(\varphi_i^w) = \begin{cases}
			1 & \text{ if } v=w,\\
			0 & \text{ otherwise.}
		\end{cases} \label{eq:sepprop}
	\]
	Then there is a finite $\Sigma\subseteq \CH(P)$ of axioms such that 
	the following are equivalent for any $B:\gF(P)\to[0,1]$.
	\begin{enumerate}[(A)]\itemsep-2pt
		\item {\bf No Dutch book}: There is no Dutch book against $B$.
		\item {\bf Convex combination}: $B$ is a convex combination 
		of the functions in $\Hom_e(\gF(P), \gA)$.
		\item {\bf Probability}: $B$ validates the axioms $\Sigma$, and 
		if $\varphi\sim\psi$, then $B(\varphi)=B(\psi)$.
	\end{enumerate}
\end{theorem}
\begin{proof}
	(A) and (B) are equivalent by Theorem \ref{pre:conv}. For the equivalence of (B) and (C), we define first $\Sigma$. Recall first that the Lindenbaum--Tarski algebra $\gL(P)$ is finite, because both $P$ and $\gA$ are finite. 
	Let $\Sigma$ be the set consisting of the following formal equalities:
	\begin{description}\itemsep-1pt
		\item[($Ax_{\psi}$)]\quad $\psi = \sum_{w}\sum_{i\in I_w}w(\psi)\alpha_i^w\varphi_{i}^w$ \quad
		for each $\psi\in \gL(P)$,
		\item[($Ax_{1}$)]\quad $\sum_w\sum_{i\in I_w}\alpha_i^w\varphi_i^w = 1$.
	\end{description}
	We show first that every cognitive evaluation validates every equality in $\Sigma$, yielding (B)$\Rightarrow$(C). For $v\in \Hom_{e}(\gF(P),\gA)$ calculation gives
	\[
		v\big(  \sum_{w}\sum_{i\in I_w}w(\psi)\alpha_i^w\varphi_{i}^w  \big) &= 
		\sum_{w}w(\psi)\sum_{i\in I_w}\alpha_i^w v(\varphi_{i}^w) = w(\psi) \\
		v\big(\sum_w\sum_{i\in I_w}\alpha_i^w\varphi_i^w \big) &= \sum_{w}\sum_{i\in I_w}\alpha_i^w v(\varphi_{i}^w) = 1.		
	\]
	As these equalities hold for every $v\in\Hom_{e}(\gF,\gA)$, they also hold with respect to every substitution of variables. Thus, cognitive evaluations validate $\Sigma$.
	Finally, for (C)$\Rightarrow$(B), assume that $B$ validates $\Sigma$. For $w\in\Hom_e(\gF,\gA)$ 
	define
	\[
		\lambda_w = \sum_{i\in I_w}\alpha_i^wB(\varphi_i^w).
	\]
	Then $\lambda_w\geq 0$, and 
	\[
		\sum_{w}\lambda_w = \sum_w\sum_{i\in I_w}\alpha_i^wB(\varphi_i^w) \overset{\dagger}{=} 1\,.
	\]
	Here, $\dagger$ holds because $B$ validates the axiom ($Ax_{1}$).
	We claim that $B(\psi) = \sum_{w}\lambda_w w(\psi)$ for every $\psi\in\gF$.
	As $B$ assigns the same value to equivalent formulas, it is enough to 
	check this equality for $\psi\in\gL$. Unfolding the right-hand side we get
	\[
		\sum_{w}\lambda_w w(\psi) = \sum_{w}\sum_{i\in I_w}\alpha_i^wB(\varphi_i^w) w(\psi)
		\overset{\star}{=} B(\psi)\,.
	\]
	Here, $\star$ holds because  $B$ validates ($Ax_{\psi}$).
\end{proof}

\noindent The infinite version of the theorem is immediate:

\begin{theorem}
	Let $P$ be any set, and suppose that for each finite $Q\subseteq P$ the 
	assumption of Theorem \ref{thm:m1} holds. 
	Then there is countable set $\Sigma\subseteq \CH(X)$ of axioms (for $X$ countable) 
	such that the following are equivalent for any $B:\gF(P)\to[0,1]$.
	\begin{enumerate}[(A)]\itemsep-2pt
		\item {\bf No Dutch book}: There is no Dutch book against $B$.
		\item {\bf Convex combination}:	For any finite $Q\subseteq P$, $B\restr \gF(Q)$
				is in the convex hull of $\Hom_e(\gF(Q), \gA)$.
		\item {\bf Probability}: $B$ validates the axioms $\Sigma$, and 
		if $\varphi\sim\psi$, then $B(\varphi)=B(\psi)$.
	\end{enumerate}
\end{theorem}
\begin{proof}
	(A)$\Leftrightarrow$(B) is Theorem \ref{pre:conv}. For the equivalence of (B) and (C)
	take $\Sigma$ to be the union of the $\Sigma(Q)$'s for finite $Q\subseteq P$ given by
	Theorem \ref{thm:m1} (for each finite cardinality it is enough to consider one $Q$ of that size).	
\end{proof}

\noindent Let us check what axioms this theorem gives in the case of classical propositional logic. 
For $P=\{p\}$, there are two cognitive evaluations: $w_0(p)=0$, and $w_1(p)=1$. The Lindenbaum--Tarski algebra has four elements (it is the $4$-element Boolean algebra):
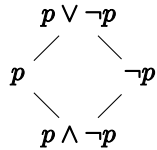
\begin{center}
\begin{tikzpicture}[scale=.4]
  \node (one) at (0,2) {$p\lor \lnot p$};
  \node (a) at (-2,0) {$p$};
  \node (d) at (2,0) {$\lnot p$};
  \node (zero) at (0,-2) {$p\land \lnot p$};
  \draw (zero) -- (a) -- (one) -- (d) -- (zero);
\end{tikzpicture}
\end{center}
The conditions of Theorem \ref{thm:m1} are satisfied with $\alpha^0=1$, $\varphi^0=\lnot p$, 
and $\alpha^1=1$, $\varphi^1=p$. Then the axioms are: 
\begin{align}
	B(p) &= B(\lnot p)w_0(p) + B(p)w_1(p) = B(p)\\
	B(\lnot p) &= B(\lnot p)w_0(\lnot p) + B(p)w_1(\lnot p) = B(\lnot p)\\
	B(p\lor \lnot p) &= B(\lnot p)w_0(p\lor\lnot p) + B(p)w_1(p\lor\lnot p) = B(\lnot p)+B(p)\\
	B(p\land\lnot p) &= B(\lnot p)w_0(p\land\lnot p) + B(p)w_1(p\land\lnot p)  = 0, \text{ and }\\
	1 &= B(p)+B(\lnot p)
\end{align}
The first four are the axioms ($Ax_{\psi}$), and the last one is ($Ax_{1}$) from Theorem \ref{thm:m1}. Clearly, the first two are redundant, and the third and last one together yields $B(\lnot p)=1-B(p)$. We would like to stress that the requirement that $B$ validates these axioms means that the equalities above hold for every substitution of any formulas in place of the variable $p$. 

For $P=\{p,q\}$ there are $4$ evaluations $w_{00}$, $w_{01}$, $w_{10}$ and $w_{11}$ (here $w_{ij}(p)=i$ and $w_{ij}(q)=j$) 
and the Lindenbaum--Tarski algebra has $16$ elements. 
The conditions of Theorem \ref{thm:m1} are satisfied with $\alpha^{ij}=1$ ($i,j<2$),
and $\varphi^{ij} = p^i\land q^j$ ($i,j<2$). We do not list all the axioms given by Theorem \ref{thm:m1} but only illustrate the method on this one example:
\[\begin{split}
	B(p) &= B(p^0\land q^0)w_{00}(p) + B(p^0\land q^1)w_{01}(p) + 
			B(p^1\land q^0)w_{10}(p) +B(p^1\land q^1)w_{11}(p) \\
		 &= B(\lnot p\land  \lnot q)\cdot 0 + B(\lnot p\land  q)\cdot 0 + 
			B(p\land \lnot q)\cdot 1 +B(p\land q)\cdot 1 \\
		 &= B(p\land \lnot q) + B(p\land q)
		 \end{split} \label{eq:axBp}
\]
The standard axiom for disjunctions (P3) can be obtained via a \emph{formal deduction} (as introduced in Subsection \ref{subsec:birkhoff}): 
\[
	B(\varphi\lor \psi) &= B((\varphi\lor \psi)\land \varphi) + 
							B((\varphi\lor \psi)\land\lnot\varphi) 
							\tag{instance of axiom \eqref{eq:axBp}} \\
						&= B((\varphi\land\varphi)\lor (\psi\land\varphi)) + 
							B((\varphi\land\lnot\varphi)\lor (\psi\land\lnot\varphi)) 
							\tag{equivalent formulas} \\
						&= B(\varphi)+B(\psi\land\lnot\varphi)
							\tag{equivalent formulas} \\
						&= B(\varphi)+B(\psi)-B(\varphi\land\psi)
						\tag{instance of ax.\eqref{eq:axBp}}
\]
For larger $P = \{p_1, \ldots, p_n\}$ the conditions of Theorem \ref{thm:m1} are still satisfied shown by
\[
	v\big(p_1^{w(p_1)}\land p_2^{w(p_2)}\land\cdots\land p_n^{w(p_n)} \big) = 
	\begin{cases}
		1 & \text{ if } v=w,\\
		0 & \text{ otherwise.}
	\end{cases}
\]
The more elements $P$ has the larger the Lindenbaum--Tarski algebra is, and the more axioms Theorem \ref{thm:m1} gives. Recall that $\CH(P)$ was the set of formal inequalities that are validated by all cognitive evaluations in $\Hom_e(\gF(P),\gA)$.
Intuitively, $\CH(P)$ is the set of all possible axioms for characterizing the convex hull of $\Hom_e(\gF(P),\gA)$. What Theorem \ref{thm:m1} provides us with is a finite set $\Sigma(P)$ of formal inequalities such that if $B:\gF(P)\to[0,1]$ validates $\Sigma(P)$, then $B$ belongs to the convex hull of $\Hom_e(\gF(P),\gA)$. But then $B$ must validate $\CH(P)$ as well. In other words, $\Sigma(P)\models \CH(P)$. 

Theorem \ref{thm:pre1} on the other hand tells us that for any $P$ the $2$-variable axioms (P1)-(P3) suffice to characterize the convex hull of the cognitive evaluations on variables $P$. In other words, there is a finite set $\Sigma\subseteq\Ineq(\{x, y\})$ such that for any finite $P$, $B:\gF(P)\to[0,1]$ validatxes $\Sigma$ if and only if $B$ belongs to the convex hull of $\Hom_e(\gF(P),\gA)$. Consequently, $\Sigma \models \CH(P)$ for all finite $P$. This is a reminiscence of compactness, and we will say that the cognitive logic has a finite base; for more discussion see the next subsection. Is it possible to get this result without a direct reference to Theorem \ref{thm:pre1}? 

In the case of propositional logic we can sketch the argument as follows. Axioms using a single variable are not sufficient because there are connectives of arity two (and if we want to replace the non-logical axiom by logical ones, then already the Boolean axioms require two variables). Assume that there is a set of axioms $\Sigma_2$ for two variables, and proceed by induction. The inductive hypothesis is that for $P$ and
$B:\gF(P)\to[0,1]$, $B\models\Sigma_2$ iff $B$ belongs to the convex hull of $\Hom_e(\gF(P),\gA)$. Let $r\notin P$ and consider $Q = P\cup\{r\}$, and $B:\gF(Q)\to[0,1]$ such that $B\models\Sigma_2$. Then every atom of the Lindenbaum--Tarski algebra $\gL(P)$ is split into two atoms in $\gL(Q)$, and we can consider $\gL(P)$ being a subalgebra of $\gL(Q)$. By the inductive hypothesis, $B$ is a probability on $\gL(P)$. If $a$ is an atom of $\gL(P)$, then axiom \eqref{eq:axBp}
applied to $a$ and $r$: $B(a) = B(a\land r)+B(a\land \lnot r)$ 
is the only constraint that $B$ needs to satisfy in order to be a probability on $\gL(Q)$. This axiom requires two variables only. 

The readers can convince themselves that analogous arguments works in the case of {\L}ukasiewicz's, Priest's, Kleene's and Symmetric logic; or just simply refer to Theorems \ref{thm:luk}, \ref{thm:sym}.

\subsection{Axioms that always work}\label{sub:always}

The axiomatization in the previous subsection depended on the assumption that homomorphisms can be separated by formulas (property \eqref{eq:sepprop}). In this subsection we show that for any given finite $P$, the set $\CH(P)$ of formal inequalities (recall: these are the formal expressions that are validated by all cognitive evaluations) always axiomatize the convex hull of the cognitive evaluations. What is more, it is always possible to chose a finite subset $\Sigma(P)\subseteq\CH(P)$ such that $\Sigma(P)$ does the job, as well. Finding elements of $\Sigma(P)$ can be done constructively/algorithmically once the values of the cognitive evaluations are known.

\begin{theorem}\label{theorem:1}
	Let $P$ be a finite set. There is a finite $\Sigma(P)\subseteq\CH(P)$ set of axioms such that for any $B:\gF(P)\to[0,1]$ the following are equivalent.
	\begin{enumerate}[(A)]\itemsep-2pt
		\item {\bf No Dutch book}: There is no Dutch book against $B$.
		\item {\bf Convex combination}: $B$ is in the convex hull of $\Hom_e(\gF(P),\gA)$.
		\item {\bf Probability}: $B$ validates the axioms $\Sigma(P)$, and 
				if $\varphi\sim\psi$, then $B(\varphi)=B(\psi)$.
	\end{enumerate}
\end{theorem}
\begin{proof}
	That (A) and (B) are equivalent is Theorem \ref{pre:conv}. We prove
	(B)$\Leftrightarrow$(C). The direction (B)$\Rightarrow$(C) is immediate 
	from that convex combinations of cognitive evaluations validate $\CH(P)$, 
	and thus $\Sigma(P)$. As for the converse direction recall that since $\gA$ 
	and $P$ are finite, there are only finitely many formulas up to 
	$\sim$-equivalence, and thus the Lindenbaum--Tarski algebra $\gL(P)$ is finite.	
	Let $\varphi_1$, $\ldots$, $\varphi_n$ be representants from each 
	$\sim$-equivalence class (elements of $\gL(P)$), and consider the $n$-dimensional 
	real valued vectors
	\[
		\vec{v} &= [v(\varphi_1), \ldots, v(\varphi_n)]\quad\text{ for } v\in V
	\]
	Let $A$ be the convex hull of $\{ \vec{v}: v\in V\}$. The convex set $A$ is
	the convex hull of finitely many points, it is a polytope, and thus it is the 
	intersection of finitely many hyperplanes. Thus, there is a finite set $R$ of
	pairs $(\vec{r}, c)$, where $\vec{r}$ is an $n$-dimensional vector, and $c\in\mathbb{R}$
	such that
	\[
		A = \bigcap_{(\vec{r},c)\in R}\big\{  \vec{x}:\; \vec{r}\cdot\vec{x}\leq c \big\}\,.
	\]
	(Here $\cdot$ is the scalar product). In other words, $\vec{x}$ belongs to 
	$A$ if and only if $\vec{r}\cdot \vec{x}\leq c$ for every pair $(\vec{r},c)\in R$.
	Every such inequality translates to the formal inequality
	\[
		r_1\cdot \varphi_1 + r_2\cdot \varphi_2 + \ldots + r_n\cdot \varphi_n \leq c\,. 
		\label{eq:ineq}
	\]
	(Here $\cdot$ is formal multiplication between the real number $r_i$, and the 
	formula $\varphi_i$). Let $\Sigma(P)$ be the finite set of these
	formal expressions. Now, each $v\in V$ validates every member of 
	$\Sigma(P)$. Without substitution this is immediate, as
	$\vec{v}$ belongs to $A$. For any proper substitution, this still holds, because
	a substituted evaluation is equal to another evaluation. So far we have that
	$\Sigma(P)\subseteq \CH(P)$.
		Let us now define 
	\[	\vec{b} &= [B(\varphi_1), \ldots, B(\varphi_n)]\,.	\]
	If $B$ is not a convex combination of the cognitive evaluations, then 
	$\vec{b}\notin A$. By the Minkowski hyperplane separation theorem there
	must be a hyperplane that separates $A$ and $\vec{b}$. In particular, 
	there must be $(\vec{r},c)\in R$ such that 
	\[
		\vec{r}\cdot\vec{v} \leq c,\quad\text{ and }\quad \vec{r}\cdot\vec{b}>c\,.
	\]
	But this means that $B$ does not validate the formal inequality \eqref{eq:ineq} 
	in $\Sigma(P)$
	corresponding to the pair $(\vec{r}, c)\in R$. And thus, $B$ does no satisfy $\Sigma(P)$.
\end{proof}

\noindent This theorem gives us a tool to determine a set of axioms: check what the polytope of the cognitive evaluations is, and find, in terms of formal inequalities, the hyperplanes that border the polytope.

The infinite version of the theorem is immediate:

\begin{theorem}\label{theorem:2}
	Let $P$ be any set. There is a countable set $\Sigma\subset\CH(X)$ of axioms 
	such that for any $B:\gF(P)\to[0,1]$ the following are equivalent.
	\begin{enumerate}[(A)]\itemsep-2pt
		\item {\bf No Dutch book}: There is no Dutch book against $B$.
		\item {\bf Convex combination}: For any finite $Q\subseteq P$, 
		$B\restr \gF(Q)$ is in the convex hull of $\Hom_e(\gF(Q), \gA)$.
		\item {\bf Probability}: $B$ validates the axioms $\Sigma$, and 
				if $\varphi\sim\psi$, then $B(\varphi)=B(\psi)$.
	\end{enumerate}
\end{theorem}
\begin{proof}
	Let $\Sigma$ be the union of the $\Sigma(Q)$'s for finite $Q\subseteq P$ given by
	Theorem \ref{thm:m1} (for each finite cardinality it is enough to consider one $Q$ of that size).	
\end{proof}

\paragraph{Getting rid of the logical axiom.} 
Theorem \ref{theorem:2} always gives an axiomatization of the convex hull of the cognitive evaluations. The only ``logical axiom'' we assumed is that if $\varphi\sim\psi$, then $B(\varphi)=B(\psi)$. The variety generated by $\gA$ can always be axiomatized by countably many equations (recall that we assumed that there are finitely many connectives only). Consider an identity $t_1=t_2$ that is used to axiomatize the variety. Such an identity is also a formal inequality, and thus can be considered as an axiom for probability. Indeed, if $t_1=t_2$ holds in $\gA$, then $B(t_1)=B(t_2)$ is a necessary condition for $B$ to belong to the convex hull of the cognitive evaluations. As already noted before, having the $\gA$-valid identities as probability axioms ensures that $B(\varphi)=B(\psi)$ holds whenever
$\varphi\sim\psi$. The moral is that we can always replace the logical axiom with at most countably many non-logical axioms. If the variety of $\gA$ is finitely based, that is, finitely many identities characterize the variety, then that amounts to replacing the logical axiom with finitely many non-logical ones. Whether or not $\gA$ is finitely based is a difficult problem in general \cite{RossWillard}. By Baker's theorem \cite[Thm.5.1]{RossWillard}, if $\gA$ is a finite algebra of finite language, and the variety generated by $\gA$ is congruence distributive, then it is finitely based. In particular, if $\gA$ has a lattice reduct, then this is the case. In the part of the literature we are contributing to, probabilities ``typically'' are considered over algebras generated by finite lattices enriched with some further operations, and thus in such ``typical'' case we always have a finite basis, and we can get rid of the logical axiom by having finitely many extra non-logical axioms instead.

\paragraph{The extended finite basis problem.}
Now that Theorem \ref{theorem:2} always gives an axiomatization of the convex hull of the cognitive evaluations, the question arises whether such an axiomatization can be finite. The role of the logical axiom is already clarified, let us turn into the question whether finitely many logical axioms, together with ``if $\varphi\sim\psi$, then $B(\varphi)=B(\psi)$'' are enough. We call the cognitive logic \emph{finitely axiomatizable}
if there is a finite set $X$ such that for any finite $P$ and $B:\gF(P)\to[0,1]$, we have that $B$ is in the convex hull of $\Hom_e(\gF(P), \gA)$ if and only if $B$ validates $\CH(X)$. In such a case, $\Sigma(X)$ provided by Theorem \ref{theorem:1} is a finite set of axioms that characterize probability mappings over the logic (for any set $P$ of propositional letters). The examples listed in Section \ref{sec:prelim} are all finitely axiomatizable cognitive logics. In the previous subsection we hinted that in the case of classical propositional logic (as well as symmetric logic, etc.) one can have a direct argument justifying that the cognitive logic is finitely axiomatizable. On might tempted to think that a finite axiomatization is always possible, and that the number of variables used in such an axiomatization might somehow depend on the size of the algebra $\gA$ and perhaps the arity of the logical connectives. In the next subsection we show a simple example that refutes these ideas.

\subsection{Meet-semilattice logic}\label{subsec:nonfin}

In this subsection, $\gA$ is the finite algebra $\< \{0,1\}, \land\>$, where $\land$ is the usual meet operation inducing the ordering $0<1$. Thus, the only logical connective we will have is the binary $\land$. To simplify notation, we think of the elements of $\gA$ as already cognitively loaded. Formally, the cognitive load function $e:\gA\to[0,1]$ gives $e(0)=0$ and $e(1)=1$. In what follows, we skip the reference to this cognitive load function. 
For a given $P$, the elements of the Lindenbaum--Tarski algebra $\gL(P)$ are 
$p_1\land\ldots p_k$ for $p_i\in P$ ($i=1\ldots k$). For any 
evaluation $w:\gF(P)\to \{0,1\}$ we have $w(p\land q) = w(p)\cdot w(q)$. This identity is \emph{not} preserved for convex combinations of evaluations (because it is not linear). 

By Theorem \ref{theorem:1} we know that for any finite $P$ there is a 
finite $\Sigma(X)\subseteq\CH(X)$ set of axioms (with $|X|=|P|$) such that $B:\gF(P)\to[0,1]$
is in the convex hull of $\Hom_e(\gF(P),\gA)$ if and only if $B$ validates the axioms $\Sigma(P)$, and $\varphi\sim\psi$ implies $B(\varphi)=B(\psi)$. We claim that there is no finite set $\Sigma(X)$ of axioms that would work for any finite $P$. This is in contrast with all the examples listed in Section \ref{sec:prelim}. First, we provide axioms that
characterize the convex hull of the evaluations. 

\begin{theorem}\label{thm:semilat1}
	Let $P=\{p_1, \ldots, p_n\}$ be a finite set, and 
	let $X=\{x_1, \ldots, x_n\}$ be a set of formula variables.
	For any $B:\gF(P)\to[0,1]$ the following are equivalent.
	\begin{enumerate}[(A)]\itemsep-2pt
		\item {\bf No Dutch book}: There is no Dutch book against $B$.
		\item {\bf Convex combination}: $B$ is in the convex hull of $\Hom_e(\gF(P),\gA)$.
		\item {\bf Probability}: If $\varphi\sim\psi$, then $B(\varphi)=B(\psi)$, and
		$B$ validates the axioms 
		\begin{enumerate}[(i)]\itemsep-1pt
			\item $x_1\land x_2\ \leq\ x_1$ (if $n\geq 2$),
			\item For every $1\leq m\leq n$, \ 
				$\sum_{k=1}^{m}(-1)^{k+1}\sum_{1\leq i_1<\ldots<i_k\leq m}
				(x_{i_1}\land\cdots \land x_{i_k})\ \leq\ 1$.
		\end{enumerate}
	\end{enumerate}
\end{theorem}
\begin{proof}
	We assume that the belief functions $B$ we consider in this proof always
	satisfy the condition that $\varphi\sim\psi$ implies $B(\varphi)=B(\psi)$. This
	amounts to treating belief functions as being defined on the Lindenbaum--Tarski
	algebras $\gL(P)$, and considering the convex hulls of the evaluations 
	$\Hom_e(\gL(P), \gA)$. To simplify presentation we will rely on the identification 
	of belief functions $B:\gF(P)\to[0,1]$ (satisfying that $\varphi\sim\psi$ 
	implies $B(\varphi)=B(\psi)$) with $B:\gL(P)\to[0,1]$. 
	As an illustration, we start with checking the simplest cases with $1$ or $2$
	propositional letters, and then we provide the proof for the general case. \\
	
	For $P_1=\{p\}$, $\gL(P_1) = \{p\}$ and there are only two evaluations: 
	$w_0(p)=0$ and $w_1(p)=1$. The convex hull of $\Hom_e(\gL(P_1), \gA)$
	is thus the set of all belief functions $B:\gL(P)\to[0,1]$. In this case axiom (i)
	is simply $x\leq 1$.\\
	
	For $P_2=\{p, q\}$ we have four evaluations $w_{ij}$ ($i,j<2$):
	\begin{center}
			\begin{tabular}{l|c|c|c}
					     & $p$ & $q$ & $p\land q$ \\ \hline
				$w_{00}$ & $0$ & $0$ & $0$ \\
				$w_{01}$ & $0$ & $1$ & $0$ \\
				$w_{10}$ & $1$ & $0$ & $0$ \\
				$w_{11}$ & $1$ & $1$ & $1$ 
			\end{tabular}
	\end{center}
	The convex hull of these evaluations (that is, the rows of the table above) is 
	the set
	\begin{align}
		\big\{ (x,y,z):\; 0\leq z,\  0\leq x-z,\  0\leq y-z, \ x+y-z\leq 1
		\big\}\,. \label{eq:semilat1}
	\end{align}
	Using this characterization and noting that $z$ corresponds to $x\land y$, and that
	the notion of validating a set of formal inequalities involves permutation of 
	variables, we have that the set
	\begin{align}
		\Sigma = \big\{ x\land y \;\leq\; x, \quad x+y- (x\land y) \;\leq\; 1 \big\}
		\label{eq:sigma1}
	\end{align}
	axiomatizes the convex hull. I.e., for any $B:\gL(P)\to[0,1]$, 
	$B\models \Sigma$ if and only if $B$ is in the convex hull of $\Hom_e(\gL(P_2),\gA)$.\\
	
	For $P_3=\{p,q,r\}$ the elements of $\gL(P_3)$ are $p$, $q$, $r$, $p\land q$, $p\land r$, 
	$q\land r$, and $p\land q\land r$. There are $8$ evaluations $w_{H}$ ($H\subseteq P_3$),
	where $w_H(p)=1$ iff $p\in H$. Every evaluation validates the formal inequality
	\begin{align}
			x+y+z-(x\land y)-(x\land z)-(y\land z)+(x\land y\land z) \;\leq\; 1\,,
	\end{align}
	which is axiom (ii) for $m=3$. Take any $B:\gF(P_3)\to[0,1]$ that we want to 
	express as a convex sum $\sum_{H\subseteq P_3}\lambda_{H}w_{H}$. 
	The only evaluation that is
	not zero on the formula $p\land q\land r$ is $w_{pqr}$. Thus $\lambda_{pqr}$ must be
	$B(p\land q\land r)$. Then, we can proceed by induction. The evaluations that are
	not zero on $p\land q$ are $w_{pq}$ and $w_{pqr}$. Thus, 
	$B(p\land q) = \lambda_{pq} + \lambda_{pqr}$, from which 
	$\lambda_{pq} = B(p\land q)-B(p\land q\land r)$. And so on. That the $\lambda$'s
	obtained this way are indeed convex coefficients follows from the axioms (i) and (ii), 
	but we skip the details for $P_3$, because we prove this below for 
	any number of propositional	letters.\\
	
	For the general case take any arbitrary non-empty finite $P = \{p_1, \ldots, p_n\}$.
	We claim first that axiom (ii) is validated by every evaluation
	$w:\gF(P)\to\{0,1\}$. For this it is enough to check that axiom (ii) holds
	with the substitutions $p_i\mapsto x_i$. Pick an evaluation $w$. 
	This $w$ maps $\ell$ variables into $1$, and $n-\ell$ variables into $0$. 
	In evaluating the left-hand side of axiom (ii), we have $\binom{\ell}{k}$ 
	choices to for a term $p_{i_1}\land\cdots\land p_{i_k}$ to be non-zero 
	(to be one). Therefore, 
	\begin{align}
		w\left( \sum_{k=1}^{n}(-1)^{k+1}\sum_{1\leq i_1<\ldots<i_k\leq n}(p_{i_1}\land\cdots
		\land p_{i_k})  \right) = \sum_{k=1}^{n}(-1)^{k+1}\binom{\ell}{k} = 
		\begin{cases}
			0&\text{ if } \ell=0,\\
			1&\text{ otherwise}
		\end{cases}
		\ \leq 1\,.
	\end{align}
	
	Next, for $\emptyset\neq H\subseteq P$ let us write $\varphi_H = \Land_{p\in H}p$, and let 
	$w_H$ be the evaluation $w_H(p) = 1$ if $p\in H$, and $w_{H}(p)=0$ if $p\notin H$. 
	For $\emptyset$, $w_{\emptyset}(p)=0$ for every $p\in P$. We need to show that
	if $B$ validates the axioms, then it belongs to the convex hull of $\Hom_e(\gF(P),\gA)$,
	that is, for every $H\subseteq P$ (including the emptyset) there is $\lambda_H$ such that
	\begin{align}
		B(\psi) = \sum_{H\subseteq P}\lambda_Hw_H(\psi)\quad\text{ for all } \psi\in \gF(P),
		\label{eq:Bkif}
	\end{align}
	and $0\leq \lambda_H\leq 1$, $\sum_{H}\lambda_H=1$. We define these $\lambda_H$'s 
	for $\emptyset\neq H\subseteq P$ as
	\begin{align}
		\lambda_H = \sum_{L\subseteq P-H}(-1)^{|L|}B(\varphi_{H\cup L})\,.
	\end{align}
	Let us check first that with this choice of $\lambda$'s \eqref{eq:Bkif} holds. 
	Take any formula $\psi$. Then there is a non-empty $H\subseteq P$ such that
	$\psi\sim\varphi_H$. The evaluations that are non-zero on $\varphi_H$ are
	exactly the evaluations $w_{H\cup K}$ for some $K$. The right-hand side of \eqref{eq:Bkif}
	is
	\begin{align}
		\sum_{A\subseteq P}\lambda_Aw_A(\varphi_H) = 
		\sum_{K\subseteq P-H}\lambda_{H\cup K}w_{H\cup K}(\varphi_H) = 
		\sum_{K\subseteq P-H}\left(\sum_{L\subseteq P-(H\cup L)} (-1)^{|L|}
		B(\varphi_{H\cup K\cup L}) \right)\,.
	\end{align}
	After expanding this latter sum, then only term which is not canceled is $B(\varphi_H)$. 
	This establishes \eqref{eq:Bkif}. 
	
	It remained to show that the $\lambda_H$'s are convex coefficients. The
	inequality $0\leq \lambda_H$ follows from axiom (i). As for $\sum_H\lambda_H=1$, 
	note that
	\begin{align}
		\sum_{\emptyset\neq H\subseteq P} \lambda_H = 
		\sum_{\emptyset\neq H\subseteq P} \left(
		\sum_{L\subseteq P-H} (-1)^{|L|}B(\varphi_{H\cup L}) \right)\ 
		\leq \ 1\,,
	\end{align}
	by axiom (ii). If this sum is strictly smaller than $1$, then choose 
	$\lambda_{\emptyset}$ such that the sum $\sum_{H\subseteq P}\lambda_H$ becomes $1$. The 
	choice of $\lambda_{\emptyset}$ does not affect the sum \eqref{eq:Bkif}, 
	because $w_{\emptyset}$ maps everything into $0$.
\end{proof}

\begin{theorem}\label{thm:semilat2}
	Let $P$ be any set, and let $X=\{x_i: i<\omega\}$ be a set of variables.
	For any $B:\gF(P)\to[0,1]$ the following are equivalent.
	\begin{enumerate}[(A)]\itemsep-2pt
		\item {\bf No Dutch book}: There is no Dutch book against $B$.
		\item {\bf Convex combination}: For any finite $Q\subseteq P$, 
		$B\restr \gF(Q)$ is in the convex hull of $\Hom_e(\gF(Q), \gA)$.
		\item {\bf Probability}: If $\varphi\sim\psi$, then $B(\varphi)=B(\psi)$, and
		$B$ validates the axioms 
		\begin{enumerate}[(i)]\itemsep-1pt
			\item $x_1\land x_2\ \leq\ x_1$,
			\item For every $n$, \ 
				$\sum_{k=1}^{n}(-1)^{k+1}\sum_{1\leq i_1<\ldots<i_k\leq n}
				(x_{i_1}\land\cdots \land x_{i_k})\ \leq\ 1$.
		\end{enumerate}
	\end{enumerate}
\end{theorem}
\begin{proof}
	Immediate from Theorem \ref{thm:semilat1}.
\end{proof}

Theorem \ref{thm:semilat1} gives for every finite $P$ a set of axioms $\Sigma(X)$ 
for $|X|=|P|$. Next, we show that there is no finite set $\Sigma(X)$ of axioms 
that would work for any finite $P$. In other words, the cognitive logic we consider here
is not finitely axiomatizable. 

\begin{theorem}
	There is no finite set $X$ and (a possible infinite) set $\Sigma(X)$ of axioms
	such that for any finite $P$ and $B:\gF(P)\to[0,1]$ the following are equivalent.
	\begin{enumerate}[(A)]\itemsep-2pt
		\item {\bf No Dutch book}: There is no Dutch book against $B$.
		\item {\bf Convex combination}: $B$ is in the convex hull of $\Hom_e(\gF(P), \gA)$.
		\item {\bf Probability}: $B$ validates the axioms $\Sigma(X)$, and 
				if $\varphi\sim\psi$, then $B(\varphi)=B(\psi)$.
	\end{enumerate}
\end{theorem}
\begin{proof}
	For $|X| = n$ let us denote the set of axioms given by Theorem \ref{thm:semilat1}
	by $\Sigma_n$. To prove the present theorem it is enough to construct for every $n\geq 2$
	a belief function $B_n:\gF(n)\to[0,1]$ such that $B_n$ validates the axioms $\Sigma_{n-1}$,
	but does not validate $\Sigma_n$. This would imply 
	that there is no finite $n$ such that the axioms $\Sigma_n$ would work 
	for any finite $P$. As each $\Sigma_n$ axiomatizes the convex hull of 
	$\Hom_e(\gF(n),\gA)$, no other set of axioms $\Sigma(X)$ could exists that would axiomatize
	the convex hulls for all $P$.

	Let us start with the case $P_3=\{p,q,r\}$ as an illustration. 
	Let $B:\gL(P_3)\to[0,1]$ be defined as
	\begin{align}
		B(p)=B(q)=B(r)=\frac{1}{2},\quad B(p\land q) = B(p\land r) = B(q\land r) = 
		B(p\land q\land r) = 0\,.
	\end{align}
	We claim first that $B$ validates $\Sigma_2$. By symmetry, 
	it is enough to check
	\begin{align}
		B(p\land q)\leq B(p),\quad B(p)+B(q)-B(p\land q) \leq 1,\quad
		B(p\land q\land r) \leq B(p\land q),	
	\end{align}
	and these inequalities indeed hold. On the other hand, $B$ does not
	validate $\Sigma_3$ and in particular 
	axiom (ii) for $n=3$ (in the statement of Theorem \ref{thm:semilat2}), because
	\begin{align}
		B(p)+B(q)+B(r)-B(p\land q)-B(p\land r)-B(q\land r)+B(p\land q\land r) = 
		3\cdot\frac{1}{2} \not\leq 1\,.
	\end{align}
	Therefore, $B$ is not in the convex hull of $\Hom_e(\gF(P_3), \gA)$. 
	Thus, the two variable inequalities $\Sigma_2$ that axiomatize the convex hull 
	of $\Hom_e(\gF(P_2), \gA)$	do not axiomatize that of $\Hom_e(\gF(P_3), \gA)$. \\
	
	In the general case let $P_n = \{p_1,\ldots, p_n\}$. 
	Let now $B:\gL(P_n)\to[0,1]$ be defined as
	\begin{align}
		B(p) = \frac{1}{n-1} \text{ for every } p\in P_n,\text{ and $B$ is $0$ 
		everywhere else}.
	\end{align}
	Then $B$ does not validate the $n$-variable instance of axiom (ii), because
	\begin{align}
		\sum_{i=1}^{n}B(p_i) - \sum_{i< j}B(p_i\land p_j) + 
		\sum_{i< j< k}B(p_i\land p_j\land p_k) - \ldots \pm B(p_1\land\cdots\land p_n)
		 = \frac{n}{n-1} \ \not\leq \ 1\,.
	\end{align}
	On the other hand, $B$ validates every at most $n-1$-variable instances of axiom (ii). 
	This is because in each such instance there is at most $n-1$ terms not involving $\land$, 
	and these are the only terms that are not zero. 
	Thus, we sum up at most $n-1$ times $\frac{1}{n-1}$, which is at most $1$.
	Consequently, $\Sigma_n$ is not derivable from $\Sigma_{n-1}$. 
\end{proof}

Finally, we show that the inequalities in the axioms cannot be replaced by equalities. This is also in sharp contrast with the examples listed in Section \ref{sec:prelim}.

\begin{theorem}
	For finite $|P|>1$, the convex hull of $\Hom_e(\gF(P), \gA)$ can not be axiomatized by
	formal equalities.
\end{theorem}
\begin{proof}
	For $|X| = n$ let us denote the set of axioms given by Theorem \ref{thm:semilat1}
	by $\Sigma_n$. These axioms describe the convex hull of $\Hom_e(\gF(n), \gA)$. 
	It is enough to show that the formal inequalities in $\Sigma_n$ cannot be expressed by
	formal equalities; and it is enough to prove the theorem for $|X|=2$.
	In particular, we show that $x+y-x\land y\leq 1$ (axiom (ii) in $\Sigma_2$)
	is not equivalent to any equality. For, the general form of a formal equality in two
	variables is $a\cdot x + b\cdot y + c\cdot (x\land y) = d$. Applying the
	evaluations to the left-hand side results in
	\begin{align}
		& w_{\emptyset}\big( a\cdot x + b\cdot y + c\cdot (x\land y)  \big) = 0 \\
		& w_{x}\big( a\cdot x + b\cdot y + c\cdot (x\land y)  \big) = a \\
		& w_{y}\big( a\cdot x + b\cdot y + c\cdot (x\land y)  \big) = b \\
		& w_{xy}\big( a\cdot x + b\cdot y + c\cdot (x\land y)  \big) = a+b+c
	\end{align}
	If such an equality	is satisfied by every evaluations, then we must have
	$d=0$, and thus $a=0$, $b=0$ and $c=0$. But then the formal equality is the
	trivial $0=0$, which is certainly not equivalent to any of the axioms. 	
\end{proof}

\end{document}

\end{document}